\documentclass[graybox]{svmult}

\usepackage{type1cm}        
\usepackage{makeidx}         
\usepackage{graphicx}        
\usepackage{multicol}        
\usepackage[bottom]{footmisc}
\usepackage{algorithm}
\usepackage{algorithmic}
\usepackage{newtxtext}       %
\usepackage{newtxmath}       
\usepackage{pgf,tikz}
\usepackage{color}
\definecolor{Martin}{rgb}{1, 0, 0}
\definecolor{Tommaso}{rgb}{0, 0, 1}
\definecolor{Roland}{rgb}{1, 0, 1}

\makeindex             
\newcommand{\Stek}{\mathcal{S}}

\begin{document}

\title*{A numerical algorithm based on probing to find optimized transmission conditions}
\author{Martin J. Gander \and Roland Masson \and Tommaso Vanzan}
\institute{Martin J. Gander \at Section de math\'{e}matiques, Universit\'{e} de Gen\`{e}ve, \email{martin.gander@unige.ch}
\and Roland Masson \at Université Côte d’Azur, CNRS, Inria, LJAD,  \email{roland.masson@unice.fr} \and Tommaso Vanzan \at CSQI Chair, Institute de math\'{e}matiques, Ecole Polytechnique F\'{e}d\'{e}rale de Lausanne, \email{tommaso.vanzan@epfl.ch}} 
%
%
\maketitle

\abstract*{Optimized Schwarz Methods (OSMs) are based on optimized
  transmission conditions along the interfaces between the
  subdomains. Optimized transmission conditions are derived at the
  theoretical level, using techniques developed in the last
  decades. The hypothesis behind these analyses are quite strong, so
  that the applicability of OSMs is still limited. In this manuscript,
  we present a numerical algorithm to obtain optimized transmission
  conditions for any given problem at hand. This algorithm requires
  few subdomain solves to be performed in an offline phase. This
  additional cost is usually negligible due to the resulting faster
  convergence, even in a single-query context.}

\section{Motivation}

Optimized Schwarz Methods (OSMs) are very versatile: they can be
used with or without overlap, converge faster compared to
other domain decomposition methods \cite{gander2006optimized}, are
among the fastest solvers for wave problems \cite{gander2019class},
and can be robust for heterogeneous problems
\cite{gander2019heterogeneous}. This is due to their general
transmission conditions, optimized for the problem at hand. Over
the last two decades such conditions have been derived for many
Partial Differential Equations (PDEs), see
\cite{gander2019heterogeneous} for a review.

Optimized transmission conditions can be obtained by diagonalizing
  the OSM iteration using a Fourier transform for two subdomains with a straight interface. This works surprisingly well, but there are important
cases where the Fourier approach fails: geometries with curved
interfaces (there are studies for specific geometries, e.g.
\cite{gigante2020optimized,gander2017optimized,Gander2014OptimizedSM}), and heterogeneous
couplings when the two coupled problems are quite different in terms
of eigenvectors of the local Steklov-Poincar\'{e} operators
\cite{gander2018derivation}. There is therefore a great need for
numerical routines which allow one to get cheaply optimized
transmission conditions, which furthermore could then lead to OSM black-box solvers. Our goal is to present one such procedure.
 
Let us consider the simple case of a two nonoverlapping subdomain
decomposition, that is $\Omega=\Omega_1\cup \Omega_2$, $\Omega_1\cap
\Omega_2=\emptyset$, $\Gamma:=\overline{\Omega_1}\cap
\overline{\Omega_2}$, and a generic second order linear PDE
\begin{equation}\label{Vanzan_mini_02_eq:PDE}
\mathcal{L}(u)=f,\text{ in }\Omega, \quad u=0\text{ on }\partial \Omega.
\end{equation}
The operator $\mathcal{L}$ could represent a homogeneous problem, i.e. the same PDE over the whole domain, or it could have
discontinuous coefficients along $\Gamma$, or even represent a
heterogeneous coupling.  Starting from two initial guesses
$u_1^0,u_2^0$, the OSM with double sided zeroth-order transmission
conditions computes at iteration $n$
\begin{equation}\label{Vanzan_mini_02_eq:OSM}
\begin{array}{r l l l l}
\mathcal{L} (u^n_1) &=0\quad & \text{on}\quad  \Omega_1,\quad
(\partial_{n_1} +s_1) u^n_1&=(\partial_{n_1} +s_1)u^{n-1}_2\quad & \text{on} \quad  \Gamma,\\
\mathcal{L} (u^n_2) &=0\quad & \text{on}\quad  \Omega_2,\quad
(\partial_{n_2} +s_2)u^n_2&=(\partial_{n_2} +s_2)u^{n-1}_1\quad & \text{on} \quad  \Gamma,
\end{array}
\end{equation}
where $s_1,s_2\in \mathbb{R}$ are the parameters to optimize.

At the discrete level, the original PDE \eqref{Vanzan_mini_02_eq:PDE} is equivalent
to the linear system
\begin{equation*}
\begin{pmatrix}
A_{II}^1 & 0 & A^1_{I\Gamma}\\
0 & A_{II}^2 &  A^2_{I\Gamma}\\
A^1_{\Gamma I}& A^2_{\Gamma I} & A_{\Gamma \Gamma}
\end{pmatrix}\begin{pmatrix}
\mathbf{u}_{1}\\ \mathbf{u}_{2}\\\mathbf{u}_\Gamma
\end{pmatrix}=\begin{pmatrix}
\mathbf{f}_1\\ \mathbf{f}_2\\\mathbf{f}_{\Gamma}
\end{pmatrix},
\end{equation*}
where the unknowns are split into those interior to domain
$\Omega_i$, that is $\mathbf{u}_i$, $i=1,2$, and those lying on
the interface $\Gamma$, i.e. $\mathbf{u}_\Gamma$. It is well known that the Dirichlet-Neumann and
Neumann-Neumann methods can be seen as Richardson type methods to
solve the discrete Steklov-Poincar\'{e} equation
\[\Sigma \mathbf{u}_\Gamma=\mathbf{\mu},\]
where $\Sigma:=\Sigma_1+\Sigma_2$,
$\Sigma_i:=A^i_{\Gamma\Gamma}-A^i_{\Gamma
  I}(A^i_{II})^{-1}A^i_{I\Gamma}$,
$\mathbf{\mu}:=\mathbf{\mu}_1+\mathbf{\mu}_2$,
$\mathbf{\mu}_i:=\mathbf{f}^i_\Gamma-A^i_{\Gamma I}(A^i_{II})^{-1}
\mathbf{f}_i$, $i=1,2$.  It is probably less known that the OSM
\eqref{Vanzan_mini_02_eq:OSM} can be interpreted as an Alternating Direction Implicit
scheme (ADI, see e.g. \cite{axelsson_1994}), for the solution of
the continuous Steklov-Poincar\'{e} equation. This interesting point
of view has been discussed in
\cite{agoshkov1990generalized,discacciati2004domain}. At the discrete
level, it results in the equivalence between a discretization of
\eqref{Vanzan_mini_02_eq:OSM} and the ADI scheme
\begin{equation*}
(s_1 E +\Sigma_1)\lambda^{n+\frac{1}{2}}=(s_1 E -\Sigma_2)\lambda^{n} +\mu,\quad
(s_2 E +\Sigma_2)\lambda^{n+1}=(s_2 E -\Sigma_1)\lambda^{n+\frac{1}{2}} +\mu,
\end{equation*}
where $E$ is either the mass matrix on $\Gamma$ using a Finite Element
discretization, or simply an identity matrix using a Finite Difference
stencil. From now on, we will replace $E$ with the identity $I$
without loss of generality.  Working on the error equation, the
iteration operator of the ADI scheme is
\begin{equation}\label{Vanzan_mini_02_eq:iteration_operator}
T(s_1,s_2):=(s_2I+ \Sigma_2)^{-1}(s_2I- \Sigma_1)(s_1I+ \Sigma_1)^{-1}(s_1I- \Sigma_2),
\end{equation}
and one would like to minimize the spectral radius,
$\min_{s_1,s_2} \rho(T(s_1,s_2))$.  It would be natural to use the
wide literature available on ADI methods to find the optimized
parameters $s_1,s_2$ for OSMs. Unfortunately, the ADI literature
contains useful results only in the case where $\Sigma_1$ and
$\Sigma_2$ commute, which is quite a strong assumption. In our
context, the commutativity holds for instance if $\Omega_1=\Omega_2$
and $\mathcal{L}$ represents a homogeneous PDE. Under these
hypotheses, Fourier analysis already provides good estimates of the
optimized parameters. Indeed it can be shown quite generally that the
Fourier analysis and ADI theory lead to the same estimates. Without
the commutativity assumption, the ADI theory relies on rough upper
bounds which do not lead to precise estimates of the optimized
parameters. For more details on the links between ADI methods and OSMs we refer to \cite[Section 2.5]{vanzan_thesis}.

Let us observe that if one used more general transmission conditions
represented by matrices $\widetilde{\Sigma}_1$ and
$\widetilde{\Sigma}_2$, \eqref{Vanzan_mini_02_eq:iteration_operator} becomes
\begin{equation*}
T(\widetilde{\Sigma}_1,\widetilde{\Sigma}_2)=(\widetilde{\Sigma}_2+ \Sigma_2)^{-1}(\widetilde{\Sigma}_2- \Sigma_1)(\widetilde{\Sigma}_1+ \Sigma_1)^{-1}(\widetilde{\Sigma}_1- \Sigma_2).
\end{equation*}
Choosing either $\widetilde{\Sigma}_1=\Sigma_2$ or
$\widetilde{\Sigma}_2=\Sigma_1$ leads to $T=0$, and thus one obtains
that the local Steklov-Poincar\'e operators are optimal transmission
operators \cite{Nataf}.

\section{An algorithm based on probing}

Our algorithm to find numerically optimized transmission
conditions has deep roots in the ADI interpretation of the OSMs
and it is based on the probing technique.  By probing, we mean the
numerical procedure through which we estimate a generic matrix $G$ by
testing it over a set of vectors. In mathematical terms, given a set
of vectors $\mathbf{x}_k$ and $\mathbf{y}_k:=G\mathbf{x}_k$, $k\in
\mathcal{K}$, we consider the problem
\begin{equation}\label{Vanzan_mini_02_eq:probing_definition}
\text{Find } \widetilde{G}\quad \text{such that}\quad \widetilde{G}\mathbf{x}_i=\mathbf{y}_i, \forall i \in \mathcal{I}.
\end{equation}
As we look for matrices with some nice properties ( diagonal,
tridiagonal, sparse...), problem \eqref{Vanzan_mini_02_eq:probing_definition} does
not always have a solution. Calling $D$ the set of admissible
matrices, we prefer to consider the problem
\begin{equation}\label{Vanzan_mini_02_eq:probing_definition_2}
\min_{\widetilde{G} \in D} \max_{k\in \mathcal{K}} \|\mathbf{y}_k-\widetilde{G}\mathbf{x}_k\|.
\end{equation}

Having remarked that the local Steklov-Poincar\'{e} operators
represent optimal transmission conditions, it would be natural to
approximate them using probing.  Unfortunately, this idea turns out to
be very inefficient.  To see this, let us carry out a continuous
analysis on an infinite strip, $\Omega_1=(-\infty,0)\times(0,1)$ and
$\Omega_2=(0,\infty)\times(0,1)$. We consider the Laplace equation
and, denoting with $\Stek_i$ the continuous Steklov-Poincar\'{e}
operators, due to symmetry we have $\Stek_1=\Stek_2=:\Stek_e$. In
this simple geometry, the eigenvectors of $\Stek_e$ are $v_k=\sin(k\pi
y)$, $k\in \mathbb{N}^+$ with eigenvalues $\mu_k=k\pi$ so that
$\Stek_e v_k= \mu_k v_k =:y_k$, see \cite{gander2006optimized}.  We
look for an operator $S=sI$, $s \in \mathbb{R}^{+}$, which
corresponds to a Robin transmission condition with parameter $s$. As
probing functions, we choose the normalized functions $v_k$,
$k=1,...,N_h$, where $N_h$ is the number of degrees of freedom on the
interface. Then \eqref{Vanzan_mini_02_eq:probing_definition_2} becomes
 \begin{equation}\label{Vanzan_mini_02_eq:sanity_check}
\begin{aligned}
\min_{S=s I,\text{ }s\in \mathbb{R}^+} \max_{k \in [1,N_h]} \|y_k-Sv_k\|
=\min_{s\in \mathbb{R}^+} \max_{k \in [1,N_h]} \|\mu_k v_k -s v_k\|= \min_{s\in \mathbb{R}^+} \max_{k \in [1,N_h]} |k\pi -s |.
\end{aligned}
\end{equation}
The solution of \eqref{Vanzan_mini_02_eq:sanity_check} is $s^*=\frac{N_h\pi}{2}$
while, according to a Fourier analysis and numerical evidence
\cite{gander2006optimized}, the optimal parameter is
$s^{\text{opt}}=\sqrt{N_h\pi}$. This discrepancy is due to the fact
that problem \eqref{Vanzan_mini_02_eq:sanity_check} aims to make the parenthesis
$(s_{i}I-\Sigma_{3-i})$, $i=1,2$ as small as possible, but it
completely neglects the other terms $(s_{i}I+\Sigma_{i})$.

This observation suggests to consider the minimization problem
\begin{equation}\label{Vanzan_mini_02_eq:opt1}
  \begin{aligned}
    \min_{\widetilde{\Sigma}_1,\widetilde{\Sigma}_2 \in D} \max_{k\in \mathcal{K}}
        \textstyle
\frac{\|\Sigma_2\mathbf{x}_k-\widetilde{\Sigma}_1\mathbf{x}_k\|}{\|\Sigma_1\mathbf{x}_k+\widetilde{\Sigma}_1\mathbf{x}_k\|}\frac{\|\Sigma_1\mathbf{x}_k-\widetilde{\Sigma}_2\mathbf{x}_k\|}{\|\Sigma_2 \mathbf{x}_k+\widetilde{\Sigma}_2\mathbf{x}_k\|}.
\end{aligned}
\end{equation}
We say that this problem is consistent in the sense that, assuming
$\Sigma_1,\Sigma_2$ share a common eigenbasis
$\left\{\mathbf{v}_k\right\}_k $ with eigenvalues
$\left\{\mu_k^i\right\}$, $\widetilde{\Sigma}_i=s_iI$, $i=1,2$,
$k=1,\dots,N_h$, then choosing $\mathbf{x}_k=\mathbf{v}_k$, we have
\begin{equation*}
\begin{aligned}
  \min_{\widetilde{\Sigma}_1,\widetilde{\Sigma}_2 \in D} \max_{k\in \mathcal{K}}
          \textstyle
          \frac{\|\Sigma_2\mathbf{x}_k-\widetilde{\Sigma}_1\mathbf{x}_k\|}{\|\Sigma_1\mathbf{x}_k+\widetilde{\Sigma}_1\mathbf{x}_k\|}\frac{\|\Sigma_1\mathbf{x}_k-\widetilde{\Sigma}_2\mathbf{x}_k\|}{\|\Sigma_2 \mathbf{x}_k+\widetilde{\Sigma}_2\mathbf{x}_k\|}=\displaystyle\min_{s_1,s_2}\max_{k \in \mathcal{I}}\in
                  \textstyle
                  \left| \frac{s_1-\mu^2_k}{s_1+\mu^1_k}\frac{s_2-\mu^1_k}{s_2+\mu^2_k}\right|
                  \displaystyle=\min_{s_1,s_2\in \mathbb{R}^+} \rho(T(s_1,s_2)),
\end{aligned}
\end{equation*}
that is, \eqref{Vanzan_mini_02_eq:opt1} is equivalent to minimize the spectral radius
of the iteration matrix.

We thus propose our numerical procedure to find optimized transmission
conditions, summarized in Steps 2-4 of Algorithm 1.
\begin{algorithm}[t]\label{Vanzan_mini_02_Alg:1}
\setlength{\columnwidth}{\linewidth}
\begin{algorithmic}[1]
\REQUIRE A set of vector $\mathbf{x}_k$, $k\in \mathcal{K}$, a characterization of $\widetilde{\Sigma}_1$, $\widetilde{\Sigma}_2$.
\STATE [Optional] For $i=1,2$, perform $N$ iterations of the power method to get approximations of selected eigenvectors $\mathbf{x}^i_k$, $i=1,2$, $k\in \mathcal{K}$. Map $\mathbf{x}^i_j$ into $\mathbf{x}_k$, for $i=1,2$, $j\in \mathcal{K}$ and $k=1,\dots,2|\mathcal{K}|$. Redefine $\mathcal{K}:=\left\{1,\dots,2 |\mathcal{K}|\right\}$.
\STATE Compute $y^i_k=\Sigma_i \mathbf{x}_k, k\in \mathcal{K}$,.
\STATE Call an optimization routine to solve \eqref{Vanzan_mini_02_eq:opt1}.
\STATE Return the matrices $\widetilde{\Sigma}_j$, $j=1,2$.
\end{algorithmic}\caption{}
\end{algorithm}
It requires as input a set of probing vectors and a characterization
for the transmission matrices $\widetilde{\Sigma}_i$, that is if the
matrices are identity times a real parameter, diagonal, or
tridiagonal, sparse etc. We then precompute the action of the local
Schur complement on the probing vectors. We finally solve
\eqref{Vanzan_mini_02_eq:opt1} using an optimization routine such as fminsearch in
\textsc{Matlab}, which is based on the Nelder-Mead algorithm.

The application of $\Sigma_i$ to a vector $\mathbf{x}_k$ requires a
subdomain solve, thus Step 2 requires $2|\mathcal{K}|$ subdomain
solves which are embarrassingly parallel. Step 3 does not require any
subdomain solves, and thus is not expensive.

As discussed in Section \ref{Vanzan_mini_02_Sec:num}, the choice of probing vectors
plays a key role to obtain good estimates. Due to the extensive
theoretical literature available, the probing vectors should be
heuristically related to the eigenvectors associated to the minimum
and maximum eigenvalues of $\Sigma_i$.  It is possible to set the
probing vectors $\mathbf{x}_k$ equal to lowest and highest Fourier
modes. This approach is efficient when the Fourier analysis itself
would provide relatively good approximations of the
parameters. However there are instances, e.g. curved interfaces
or heterogeneous problems, where it is preferable to have
problem-dependent probing vectors. We thus include an additional
optional step (Step 1), in which, starting from a given set of probing
vectors, e.g Fourier modes, we perform $N$ iterations of the power
method, which essentially correspond to $N$ iterations of the OSM, to
get more suitable problem-dependent probing vectors.  To compute the
eigenvector associated to the minimum eigenvalue of $\Sigma_i$, we
rely on the inverse power method which requires to solve a Neumann
boundary value problem.  Including Step 1, Algorithm 1 requires in
total $2|\mathcal{K}|(N+2)$ subdomain solves, where $|\mathcal{K}|$ is
the number of probing vectors in the input.

\section{Numerical experiments}\label{Vanzan_mini_02_Sec:num}

\begin{figure}[t]
\centering
\includegraphics[scale=0.3]{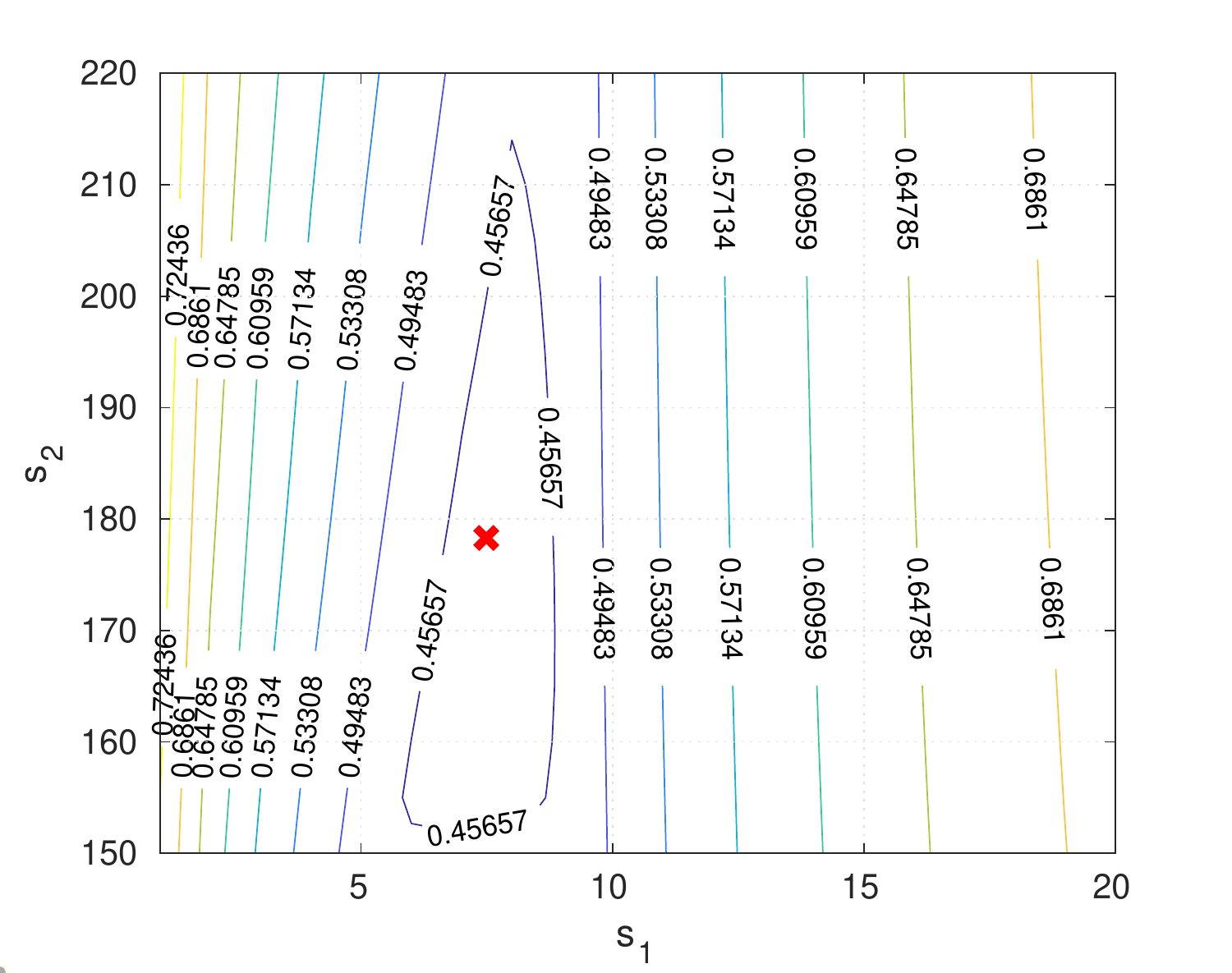}\quad
\includegraphics[scale=0.3]{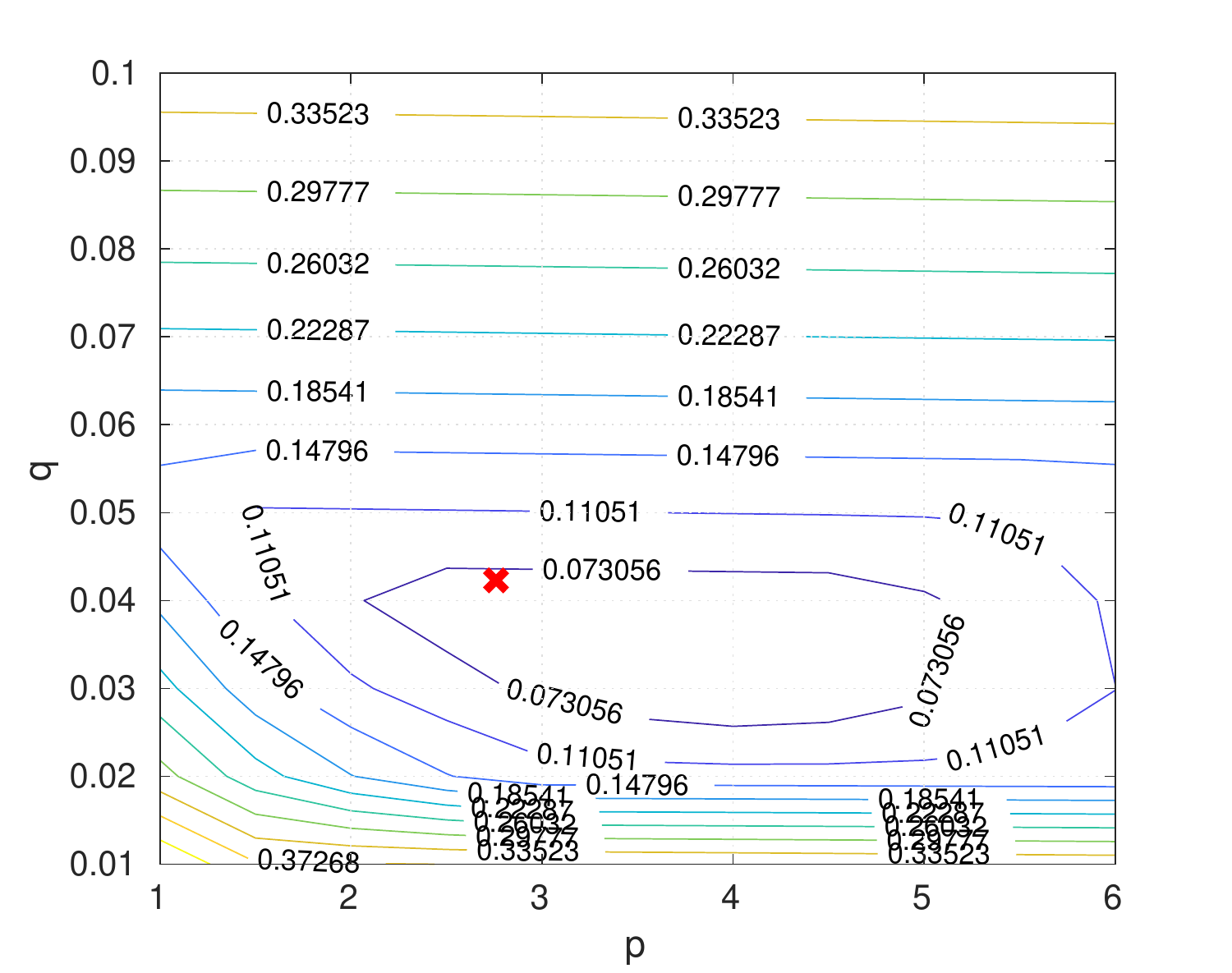}
\caption{Contour plot of the spectral radius of the iteration matrix
  $T(\widetilde{\Sigma}_1,\widetilde{\Sigma}_2)$ with
  $\widetilde{\Sigma}_i=s_i I$ (left) and of
  $T(\widehat{\Sigma}_1,\widehat{\Sigma}_2)$ with $\widehat{\Sigma}_i=
  p I +qH$ (right). The red crosses are the parameters obtained
  through Alg. 1.}\label{Vanzan_mini_02_Fig:Probing_Laplace}
\end{figure}

We start with a sanity check considering a Laplace equation on a
rectangle $\Omega$, with $\Omega_1=(-1,0)\times (0,1)$,
$\Omega_2=(0,1)\times (0,1)$ and $\Gamma=\{0\}\times (0,1)$. Given a
discretization of the interface $\Gamma$ with $N_h$ points, we choose
as probing vectors the discretization of
\begin{equation}\label{Vanzan_mini_02_eq:Probing_vectors}
x_1=\sin(\pi y),\quad x_2=\sin(\sqrt{N_h}\pi y),\quad x_3=\sin(N_h\pi y),
\end{equation}
motivated by the theoretical analysis in \cite{gander2006optimized},
which shows that the optimized parameters $s_i$ satisfy
equioscillation between the minimum, the maximum and a medium
frequency which scales as $\sqrt{N_h}$. We first look for matrices
$\widetilde{\Sigma}_i= s_iI$ representing zeroth order double sided
optimized transmission conditions. Then, we look for matrices
$\widehat{\Sigma}_i= p I +qH$, where $H$ is a tridiagonal matrix
$H:=\text{diag}(\frac{2}{h^2})-\text{diag}(\frac{1}{h^2},-1)-\text{diag}(\frac{1}{h^2},+1)$,
where $h$ is the mesh size. At the continuous level,
$\widehat{\Sigma}_i$ represent second order transmission conditions.
Fig. \ref{Vanzan_mini_02_Fig:Probing_Laplace} shows that Alg. 1 permits to obtain
excellent estimates in both cases with just three probing vectors.
We emphasize that Alg. 1 requires 6 subdomain solves, which can
be done in parallel, and leads to a convergence factor of order
$\approx 0.07$ for second order transmission conditions. It is clear
that, depending on the problem at hand, this addition of 6 subdomain
solves is negligible, considering the advantage of having such a small
convergence factor.

\begin{figure}[t]
\centering
\mbox{\begin{tikzpicture}[scale=0.5]
\draw[white] (-3.1,-2) -- (-3.1,0);
\draw[black] (-3,0) -- (0,0);
\draw[black] (-3,3) -- (0,3);
\draw[black] (-3,0) -- (-3,3);
\draw[dashed] (0,0) -- (0,3);
\node at (-1.5,1.5) {\small{$\Omega_1$}};
\node at (1.5,1.5) {\small{$\Omega_2$}};
\node at (-0.8,2.6) {\small{$\Gamma$}};
\node at (-0.5,0.4) {\small{$r$}};
\draw[<->] (0,0.4) -- (0.5,0.4);
\draw[domain=0:3,smooth,variable=\y] plot ({0.5 * (sin(deg(4/3*pi*\y)))},{\y});
\draw[domain=pi/2:-pi/2,smooth,variable=\y] plot ({3 * (cos(deg(\y)))},{1.5 * (sin(deg(\y)))+1.5});
\node at (-1,-2) {\small{$\Omega_1$}};
\draw[domain=0:1,smooth,variable=\t] plot ({2*cos(deg(2*pi*\t))-1},{1*sin(deg(2*pi*\t))-2});
\draw[black] (-1,-1) -- (2.5,-1);
\draw[black] (2.5,-1) -- (2.5,-3);
\draw[black] (2.5,-3) -- (-1,-3);
\node at (2,-2) {\small{$\Omega_2$}};
\node at (1.3,-1.5) {\small{$\Gamma$}};

\end{tikzpicture}
{\includegraphics[scale=0.3]{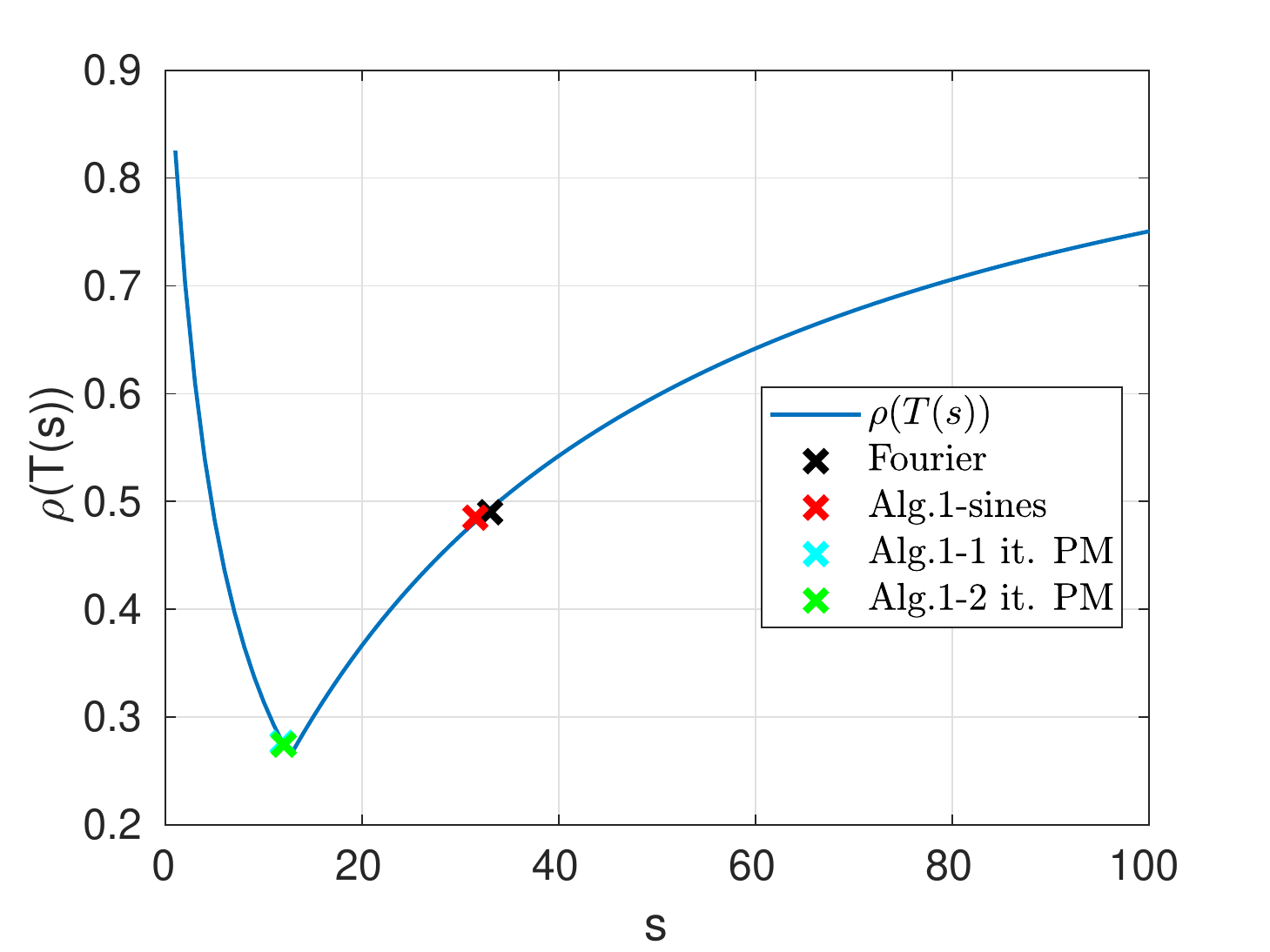}}
{\includegraphics[scale=0.3]{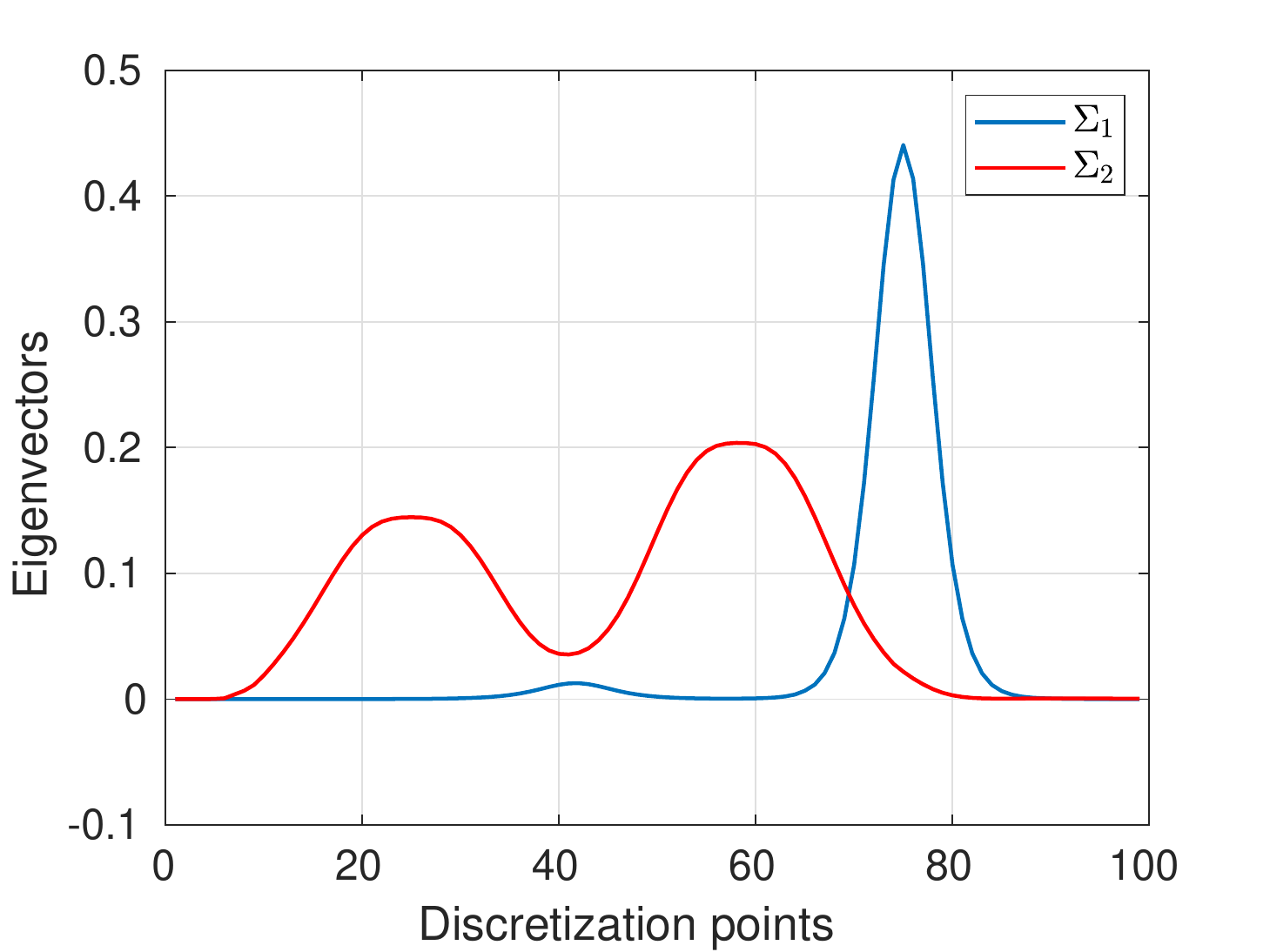}}}
\caption{Left: $\Omega$ decomposed into $\Omega_1$ and
    $\Omega_2$. Middle: optimized parameters obtained using Fourier
  analysis or Algorithm 1 with different sets of probing
  vectors. Right: eigenvectors associated to the smallest
  eigenvalues of $\Sigma_j$, $j=1,2$.}\label{Vanzan_mini_02_Fig:Probing_curved}
\end{figure}
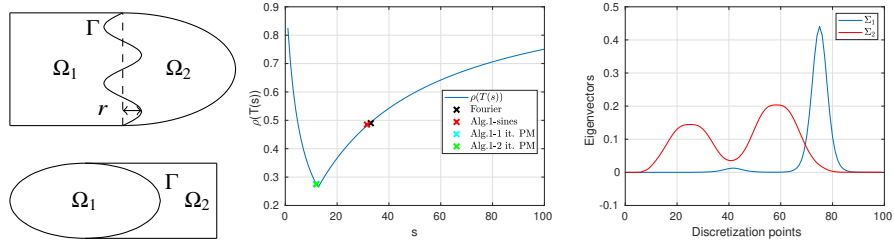

We now look at a more challenging problem. We solve a second order PDE
\begin{equation}\label{Vanzan_mini_02_eq:Probing_model}
-\nabla\cdot \nu(\mathbf{x})\nabla u+ \mathbf{a}(\mathbf{x})^\top\cdot \nabla u + \eta(\mathbf{x}) u=f\quad \text{in } \Omega,
\end{equation}
where $\Omega$ is represented in Fig. \ref{Vanzan_mini_02_Fig:Probing_curved} on the top-left.  

The interface $\Gamma$ is the parametric curve
$\gamma(t):[0,1]\rightarrow (r\sin(\widehat{k}\pi t),t)$, with $r\in
\mathbb{R}^+$. The coefficients are set to $\nu(\mathbf{x})=1$,
$\mathbf{a}(\mathbf{x})=(10(y+x^2),0)^\top$,
$\eta(\mathbf{x})=0.1(x^2+y^2)$ in $\Omega_1$, $\nu(\mathbf{x})=100$,
$\mathbf{a}(\mathbf{x})=(10(1-x),x)^\top$, $\eta(\mathbf{x})=0$ in
$\Omega_2$, $f(\mathbf{x})=x^2+y^2$ in $\Omega$. The geometric
parameters are $r=0.4$, $\widehat{k}=6$ and the interface is
discretized with $N_h=100$ points.  Driven by the theoretical analysis
\cite{gander2019heterogeneous}, we rescale the transmission conditions
according to the physical parameters, setting $S_i:=f_i(s)I$, where
$f_i:=\nu_i(s^2+ \frac{a_{i1}^2}{4\nu_i^2}
+\frac{a_{i2}^2}{4\nu_i^2} +\frac{\eta_i}{\nu_i})^{1/2}
-\frac{a_{i1}}{2}$.  The left panel of Fig. \ref{Vanzan_mini_02_Fig:Probing_curved}
shows a comparison of the optimized parameters obtained by a Fourier
analysis to the one obtained by Alg. 1 using as probing vectors
the sine frequencies \eqref{Vanzan_mini_02_eq:Probing_vectors}. It is evident that
both do not deliver efficient estimates. The failure of Alg. 1 is
due to the fact that, in contrast to the Laplace case, the sine
frequencies do not contain information about the slowest modes. On the
right panel of Fig \ref{Vanzan_mini_02_Fig:Probing_curved}, we plot the lowest
eigenvectors of $\Sigma_i$, which clearly differ significantly from
the simple lowest sine frequency.  We therefore consider Alg. 1
with the optional Step 1 and as starting probing vectors we only use
the lowest and highest sine frequencies.  The left panel of
Fig. \ref{Vanzan_mini_02_Fig:Probing_curved} shows that Alg. 1 delivers
efficient estimates with just one iteration of the power method.
Let us now study the computational cost. To solve
\eqref{Vanzan_mini_02_eq:Probing_model} up to a tolerance of $10^{-8}$ on the error,
an OSM using the Fourier estimated parameters (black cross in Fig
\ref{Vanzan_mini_02_Fig:Probing_curved}) requires 21 iterations, while only 12 are
needed by Algorithm 1 with only one iteration of the power method. In
the offline phase of Algorithm 1, we need to solve $4$ subdomain
problems in parallel in Step 1, and further $8$ subdomain problems
again in parallel in Step 2. Therefore the cost of the offline phase
is equivalent to two iterations of the OSM in a parallel
implementation, and consequently Alg. 1 is computationally
attractive even in a single-query context.

Fourier estimates depend on the choice of $k_{\min}$ and $k_{\max}$
  and in Fig. \ref{Vanzan_mini_02_Fig:Probing_curved}, we set $k_{\min}=\pi$ and
  $k_{\max}=\pi/h$. Inspired by \cite{Gander2014OptimizedSM} and a
  reviewer's comment, we optimized with $k_{\min}=\pi/|\Gamma|\approx
  \pi/4.96$ obtaining $s=14.41$, which is very close to the optimal
  $s^*$. However, rescaling $k_{\min}$ with $|\Gamma|$ is not
  generally a valid approach. Considering $\Omega_1$ as the ellipse of
  boundary $(\cos(2\pi t),0.5\sin(2\pi t))$, $t\in (0,1)$, and
  $\Omega_2=[0,2]\times [0,1]\setminus \Omega_1$, see Fig. \ref{Vanzan_mini_02_Fig:Probing_curved} bottom-left, then $s^*=40$, while
  $s_{k_{\min}=\pi}= 31.5$ and
  $s_{k_{\min}=\pi/|\Gamma|}=20.44$. Thus, rescaling $k_{\min}$
  worsens the Fourier estimate.

Next, we consider the Stokes-Darcy system in $\Omega$, with
$\Omega_1=(-1,0)\times (0,1)$, $\Omega_2=(0,1)\times (0,1)$ and
$\Gamma=\{0\}\times (0,1)$ with homogeneous Dirichlet boundary
conditions along
$\partial\Omega$. Refs. \cite{gander2018derivation,vanzan_thesis} show
that the Fourier analysis fails to provide optimized parameters since
the two subproblems do not share a common separation of variable
expansion in bounded domains, unless periodic boundary conditions are
enforced, see also \cite{gander2019heterogeneous}[Section 3.3]. Thus,
the sine functions do not diagonalize the OSM iteration, even in the
simplified domain $\Omega$ with straight interface.  Nevertheless, we
apply Alg. 1 using two different sets of sines as probing
vectors, corresponding to frequencies
$\mathcal{K}_1=\left\{1,\sqrt{N_h},N_h\right\}$ and
$\mathcal{K}_2=\left\{1,2,\sqrt{N_h},N_h\right\}$. In $\mathcal{K}_2$
the first even frequency is included because in
Ref. \cite{gander2018derivation} it was observed that the first
odd Fourier frequency converges extremely fast.

\begin{figure}[t]
\centering
\includegraphics[scale=0.3]{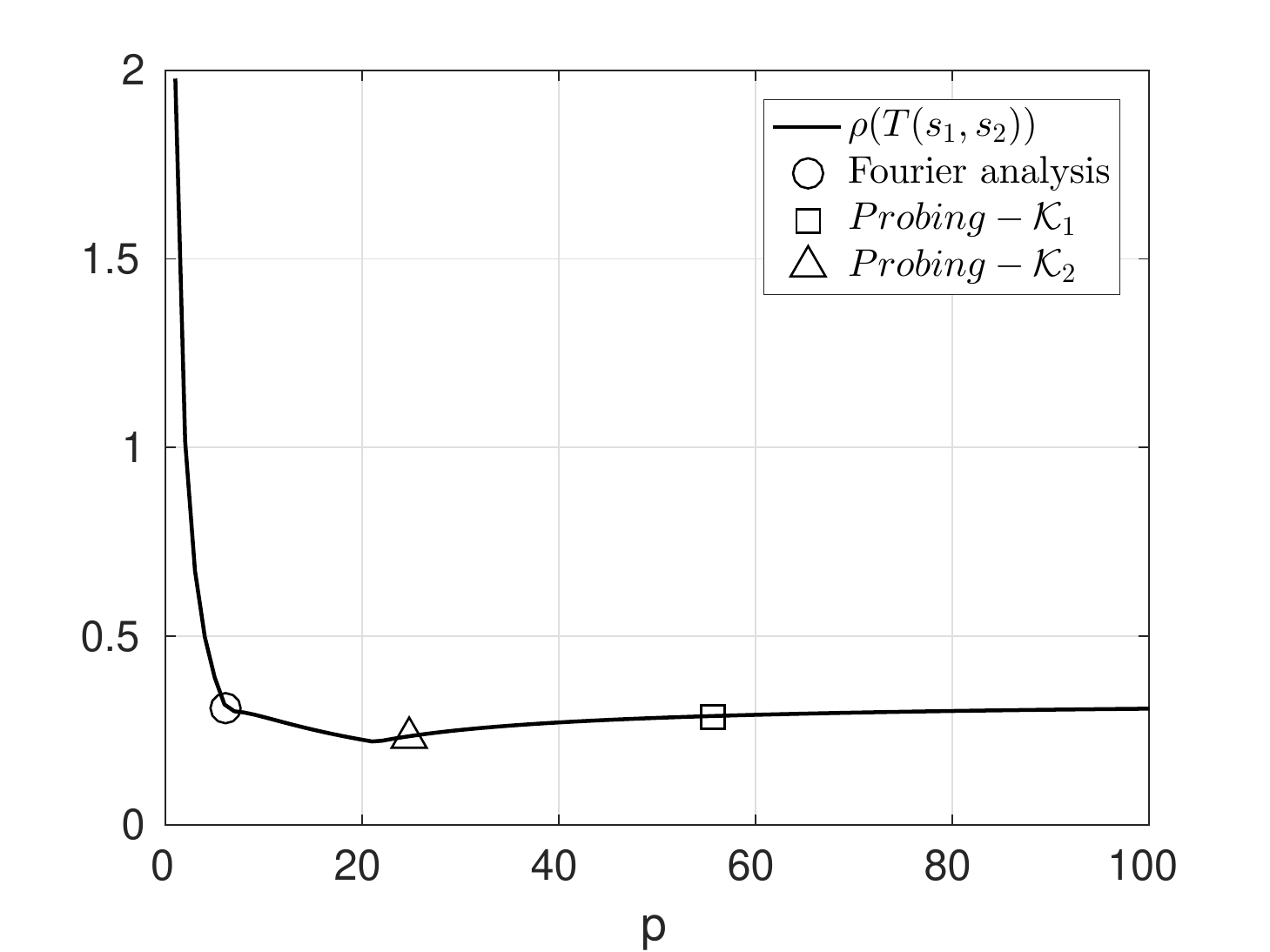}
\includegraphics[scale=0.3]{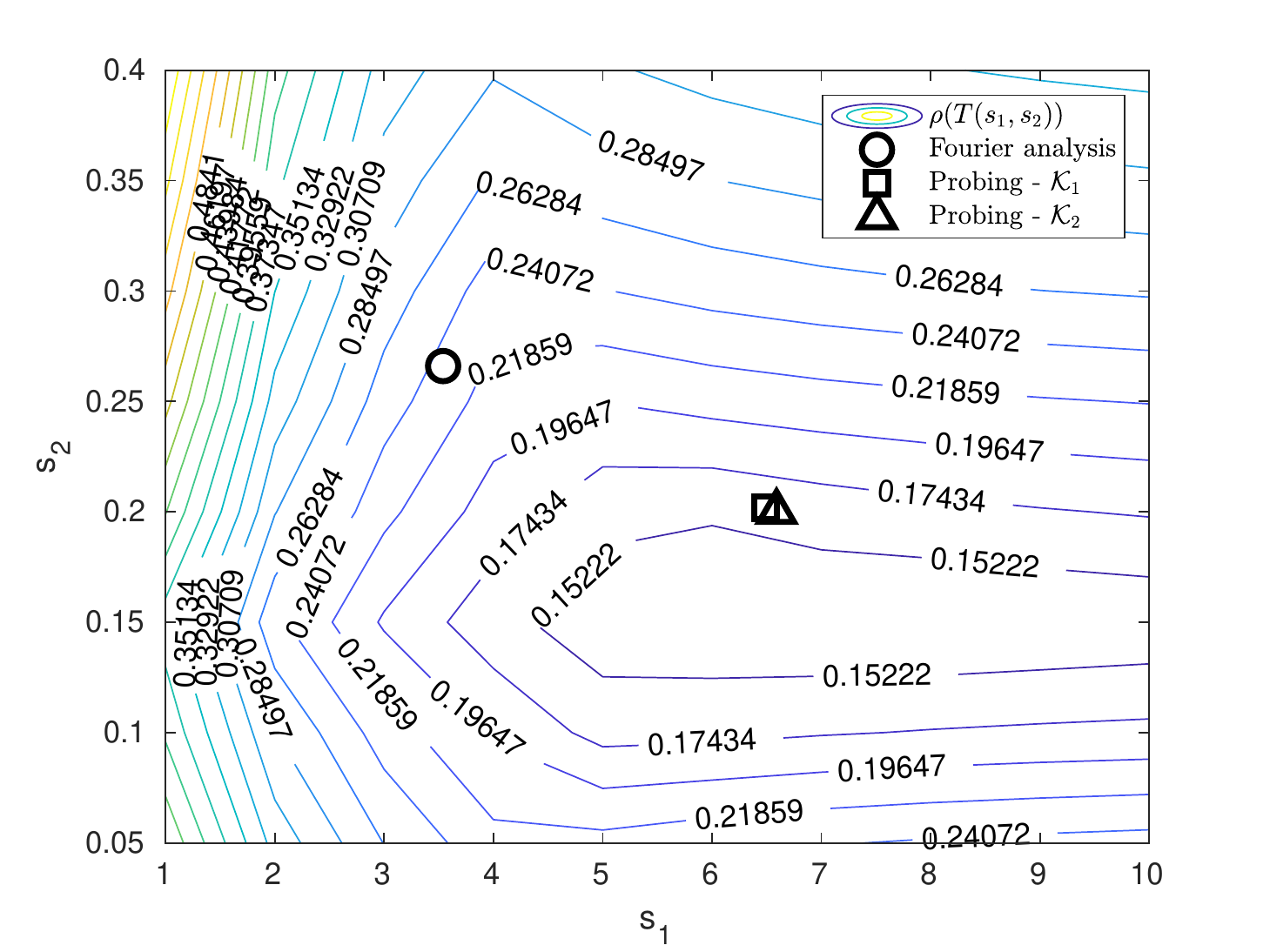}
\caption{Comparison between the optimized parameters obtained through
  Fourier analysis and Alg. 1 for single sided Robin boundary
  conditions (left) and double sided Robin boundary conditions
  (right).}\label{Vanzan_mini_02_Fig:Stokes}
\end{figure}

Fig \ref{Vanzan_mini_02_Fig:Stokes} shows the estimated parameters for single and
double sided zeroth order transmission conditions obtained through a
Fourier analysis \cite{discacciati2018optimized} and using Alg.
1.  The left panel confirms the intuition of
\cite{gander2018derivation}, that is, the first even frequency plays a
key role in the convergence.  The right panel shows that Alg. 1,
either with $\mathcal{K}_1$ or $\mathcal{K}_2$ provides better
optimized parameters than the Fourier approach.

\begin{figure}[t]
\centering
\includegraphics[scale=0.19]{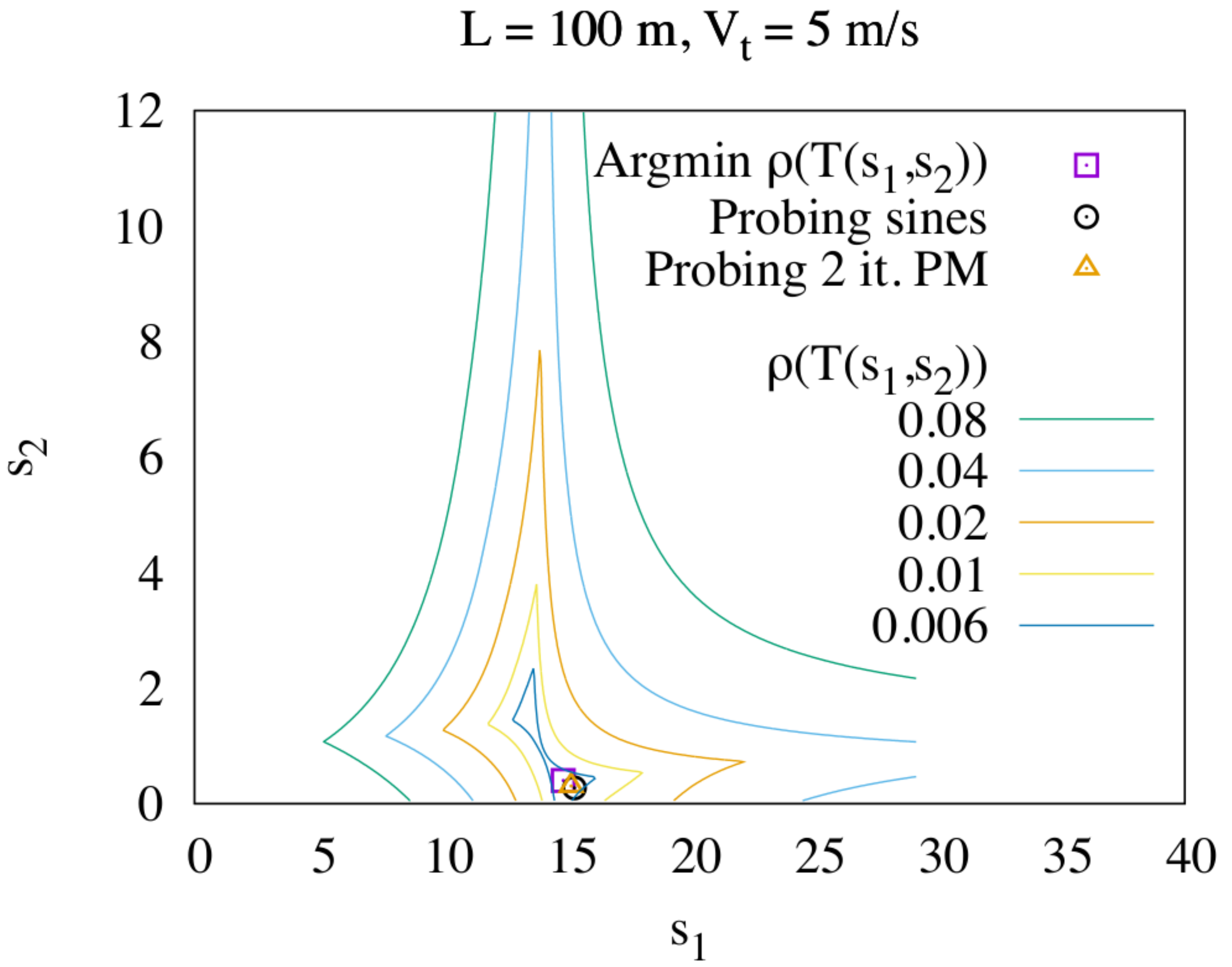}
\includegraphics[scale=0.19]{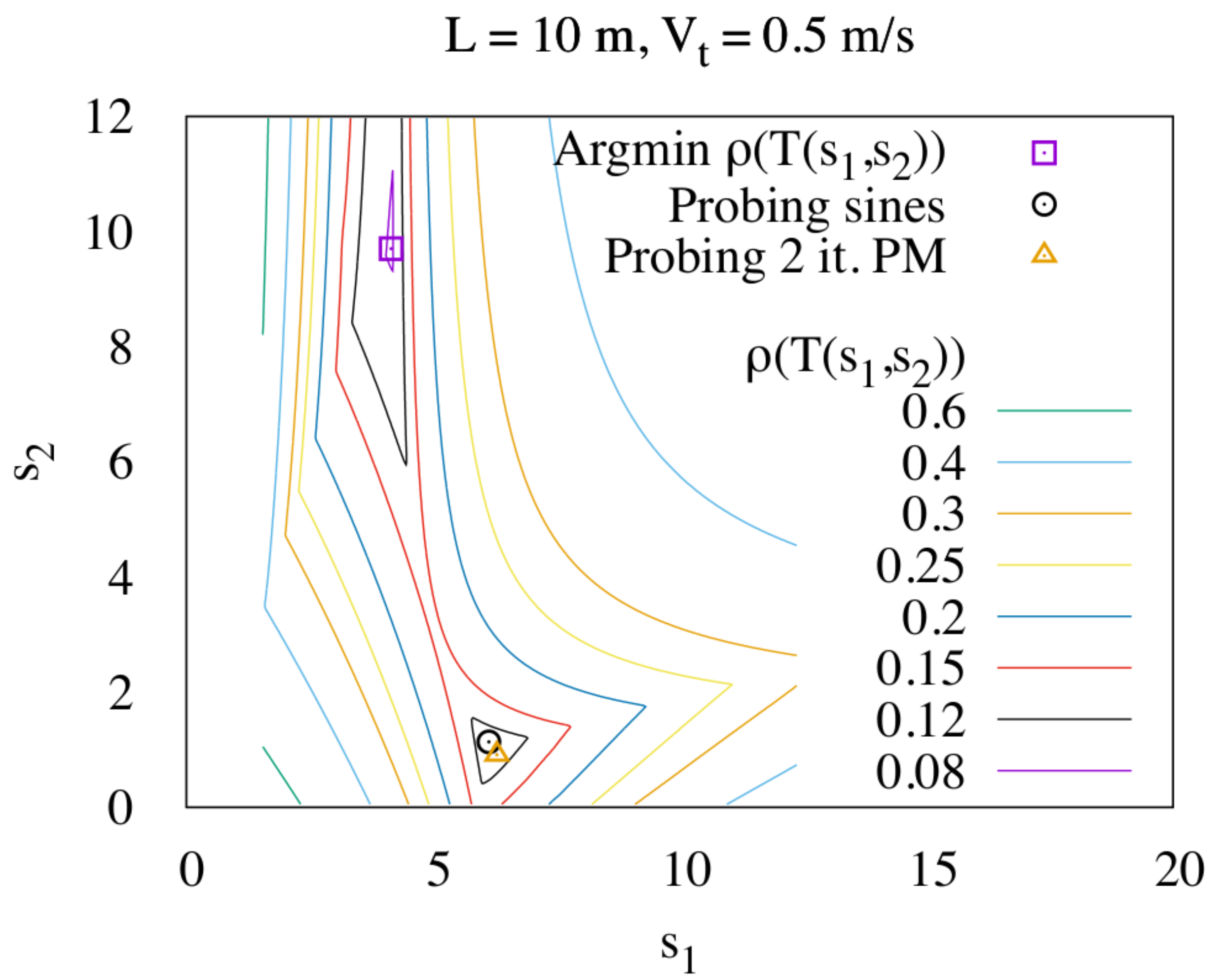}
\caption{For $L=100$ m, $\overline V_t=5$ m/s (left) and $L=10$ m,
  $\overline V_t=0.5$ m/s (right), comparison of the double sided
  Robin parameters $s_1$ and $s_2$ obtained from the probing algorithm
  using either the 3 sine vectors or the 6 vectors obtained from the
  3 sines vectors by 2 PM iterations on both sides. It is compared
  with the minimizer of the spectral radius
  $\rho(T(s_1,s_2))$. }\label{Vanzan_mini_02_Fig:pmff}
\end{figure}

Next, we consider the stationary heat transfer model coupling the
diffusion equation $\nabla\cdot \left(-\lambda \nabla
u_1(\mathbf{x})\right)=0$ in the porous medium domain
$\Omega_1=(0,L)\times (5,15)$ with the convection diffusion equation $
\nabla\cdot \left( u_2(\mathbf{x}) {\bf V}_t(y) -\lambda_t(y) \nabla
u_2(\mathbf{x}) \right)=0$ in the free flow domain
$\Omega_2=(0,L)\times (0,5)$. Both the turbulent velocity ${\bf V}_t =
(V_t(y),0)^T$ and the thermal conductivity $\lambda_t(y)$ exhibit a boundary layer at the interface
$\Gamma=(0,L)\times\{5\}$ and are computed from the Dittus-Boelter
turbulent model. Dirichlet boundary conditions are prescribed at the
top of $\Omega_1$ and on the left of $\Omega_2$, homogeneous
Neumann boundary conditions are set on the left and right of
$\Omega_1$ and at the bottom of $\Omega_2$, and a zero Fourier flux is
imposed on the right of $\Omega_2$. Flux and temperature
continuity is imposed at the interface $\Gamma$. The model is
discretized by a Finite Volume scheme on a Cartesian mesh of size
$50\times 143$ refined on both sides of the interface.  Figure
\ref{Vanzan_mini_02_Fig:pmff}
shows that the probing algorithm provides a very good approximation of
the optimal solution for the case $L=100$ m, $\overline V_t = 5$ m/s
(mean velocity) both with the 3 sine vectors
\eqref{Vanzan_mini_02_eq:Probing_vectors} and with the 6 vectors obtained from the
power method starting from the sine vectors. In the case $L=10$
m, $\overline V_t = 0.5$ m/s, the spectral radius has a narrow
valley with two minima. In that case the probing algorithm fails to
find the best local minimum but still provides a very efficient
approximation.

\end{document}